
\input amstex.tex
\documentstyle{amsppt}

\nopagenumbers

\def\DJ{\leavevmode\setbox0=\hbox{D}\kern0pt\rlap
 {\kern.04em\raise.188\ht0\hbox{-}}D}
\def\txt#1{{\textstyle{#1}}}
\def\D{\Delta} \def\e{\varepsilon} 

\def\d{{\,\roman d}}
\def\kk{\kappa_j}
  
\def\z{\zeta} \def\zx{\zeta({\textstyle{1\over2}} + ix)}
\def\zt{\zeta({\textstyle{1\over2}} + it)}
 \def\k{\kappa} \def\a{\alpha}
\def\H{H_j^3({\txt{1\over2}})} 
\font\teneufm=eufm10
\font\seveneufm=eufm7
\font\fiveeufm=eufm5
\newfam\eufmfam
\textfont\eufmfam=\teneufm
\scriptfont\eufmfam=\seveneufm
\scriptscriptfont\eufmfam=\fiveeufm
\def\mathfrak#1{{\fam\eufmfam\relax#1}}
\def\hf{{\textstyle{1\over2}}}

\font\rr=cmr9
\font\tenmsb=msbm10
\font\sevenmsb=msbm7
\font\fivemsb=msbm5
\newfam\msbfam
      \textfont\msbfam=\tenmsb
      \scriptfont\msbfam=\sevenmsb
      \scriptscriptfont\msbfam=\fivemsb
\def\Bbb#1{{\fam\msbfam #1}}

\font\fontline=cmr8
\def \NN {\Bbb N}

\def \RR {\Bbb R}

\def\rightheadline{{\hfil{\rr
On some  mean value results involving}  $ |\zeta({1\over2}+it)|$
\hfil\sevenrm\folio}}

  \def\leftheadline{{\fontline\folio\hfil{\rr
  A. Ivi\'c  }\hfil}}
  \def\emptyheadline{\hfil}
  \headline{\ifnum\pageno=1 \emptyheadline\else
  \ifodd\pageno \rightheadline \else \leftheadline\fi\fi}

\topmatter  
\title 
On some  mean value results involving $|\zt|$
\endtitle
\author    A. IVI\'C \endauthor
\dedicatory
\enddedicatory
\address{
Aleksandar Ivi\'c, Katedra Matematike RGF-a
Universiteta u Beogradu, \DJ u\v sina 7, 11000 Beograd,
Serbia (Yugoslavia).}\endaddress
\email
{\tt aivic\@rgf.bg.ac.yu, aivic\@matf.bg.ac.yu}
\endemail
\abstract
Several problems involving $E(T)$ and $E_2(T)$, the error terms 
in the mean square  and mean fourth moment formula
for $|\zt|$, are discussed. In particular it is proved that
$$
\int_0^T E(t)E_2(T)\d t \ll T^{7/4}(\log T)^{7/2}\log\log T.
$$
\endabstract
\endtopmatter

\head 
1. Introduction and statement of results
\endhead

\bigskip
Let, as usual ($\gamma = 0.5772157...$ is Euler's constant),
$$
E(T) \;=\;\int_0^T|\zt|^2\d t - T\left(\log\left({T\over2\pi}\right) 
+ 2\gamma - 1\right)
$$
denote  the error term in the mean
square formula for $|\zt|$,
and let
$$
E_2(T) \;=\; \int_0^T|\zt|^4\d t - TP_4(\log T)
$$
denote the error term in the asymptotic formula for the
fourth moment of $|\zt|$. Here $P_4(x)$ is a polynomial
of degree four in $x$ with leading coefficient $1/(2\pi^2)$ (see
[5] for the explicit evaluation of all the coefficients).
Both of these functions play an important r\^ole in
the theory of the Riemann zeta-function $\z(s)$, and 
the aim of this note is to discuss several problems
involving their mean values. Especially interesting seems the 
evaluation of the integral
$$
\int_0^T E(t)E_2(t) \d t,
$$
or (which is technically more convenient)
$$
I(T)  \;:=\; \int_T^{2T} E(t)E_2(t) \d t,\eqno(1.1)
$$
since this integral exhibits the superpositions of oscillations of
the functions $E(t)$ and $E_2(t)$. Namely both functions
are oscillating, and we have $E(t) = \Omega_\pm(t^{1/4})$ (see [1], [4])
and $E_2(t) = \Omega_\pm(t^{1/2})$ (see [4], [6] and [12]). As usual,
$f = \Omega_\pm(g)$ means that $\limsup f/g > 0$ and $\liminf f/g < 0$.
We also have (see [2] and [4])
$$
\int_0^T E^2(t) \d t = DT^{3/2} + O(T\log^4T)\eqno(1.2)
$$
with $D = 2(2\pi)^{-1/2}\z^4(3/2)/(3\z(3))$, and (see [7], [12])
$$
\int_T^{2T} E^2_2(t) \d t   \ll T^2\log^{22}T.\eqno(1.3)
$$
As usual, $f \ll g$ (same as $f = O(g)$) means that $|f(x)| \le Cg(x)$
for some $C>0$ and $x\ge x_0$).
Hence by (1.2), (1.3), and the Cauchy-Schwarz inequality for integrals 
one obtains
$$
I(T) \;\ll\; T^{7/4}\log^{11}T,\eqno(1.4)
$$
and one naturally asks whether (1.4) can be improved. This is indeed so,
as shown by the following

\bigskip
THEOREM 1. {\it We have}
$$
I(T) \;= \; \int_T^{2T} E(t)E_2(t)\d t  \;\ll\; 
T^{7/4}(\log T)^{7/2}\log\log T.\eqno(1.5)
$$
\bigskip
We remark here that an analogous formula to (1.5) holds if $E(t)$
is replaced ($d(n)$ is the number of divisors of $n$) by
$$
\D(x) = \sum_{n\le x}d(n) - x(\log x + 2\gamma -1),
$$
the error term in the classical divisor problem. 
Namely the analogues of (2.4) and (2.5) will hold for
$\D(x)$ by the Voronoi explicit formula for $\D(x)$
(see [2, Chapter 3]). Thus, following the proof of
Theorem 1, we shall obtain
$$
\int_T^{2T} \D(t)E_2(t)\d t  \;\ll\; T^{7/4}(\log T)^{7/2}\log\log T.
$$

In the course of the proof of Theorem 1 we shall encounter the function 
$$
g(t) := \hf \left({2\over\pi}\right)^{3/4}\sum_{n=1}^\infty
(-1)^nd(n)n^{-5/4}\sin(\sqrt{8\pi nt} - {\textstyle{1\over4}}\pi),\eqno(1.6)
$$
which appears in the relation
$$
\int_T^{2T} E(t)E_2(t)\d t  = 
- \int_T^{2T}t^{3/4}g(t)|\zt|^4 \d t + O(T^{3/2}\log^{10}T).\eqno(1.7)
$$
It is thus seen that the natural question of the true order of magnitude
of $I(T)$ involves the evaluation of the integral on the right-hand side
of (1.7), which contains the oscillatory function $g(t)$ (it is both
$O(1)$ and $\Omega_\pm(1)$; see [1] and [4]). Although it appears to the
author that the true order of $I(T)$ is $T^{3/2+o(1)}$, this is certainly
hard to prove. However if we set
$$
g_+(t) = \max(g(t),\,0),\qquad g_-(t) = \min(g(t),\,0),\eqno(1.8)
$$
then we have

\bigskip
THEOREM 2. {\it If $g_\pm(t)$ is given by} (1.8), {\it then for} $k = 1,2$
{\it  we have}
$$\eqalign{&
\int_T^{2T}t^{3/4}g_+(t)|\zt|^{2k}\d t  \;\ll\; T^{7/4}(\log T)^{k^2},
\cr& \int_T^{2T}t^{3/4}g_+(t)|\zt|^{2k}\d t \;\gg\; T^{7/4}(\log T)^{k^2},
\cr}\eqno(1.9)
$$
{\it and likewise}
$$\eqalign{&
\int_T^{2T}t^{3/4}g_-(t)|\zt|^{2k}\d t  \;\ll\; T^{7/4}(\log T)^{k^2},\cr&
\int_T^{2T}t^{3/4}g_-(t)|\zt|^{2k}\d t  \;\gg\; T^{7/4}(\log T)^{k^2}.\cr}
\eqno(1.10)
$$

\bigskip
Comparing Theorem 1 and Theorem 2 we see that there must be some
cancellation when we deal with $g(t) = g_+(t) + g_-(t)$ instead of only
$g_+(t)$ or $g_-(t)$. We note that,
similarly to (1.7), we obtain (with suitable $B>0$, 
see the remark after (2.2))
$$
 \int_T^{2T}t^{3/4}g(t)|\zt|^2\d t = - BT^{3/2} + O(T^{5/4}\log T).
\eqno(1.11)
$$ 
 The integrals in (1.7) and (1.11) containing the function $g(t)$ are 
similar, which is why there is reason to think that $I(T)$ is also
of the the order $T^{3/2+o(1)}$. The original motivation for the
study of $I(T)$ was to try to obtain a lower bound for the integral
on the right-hand side of (1.7). This would in turn, by the Cauchy-Schwarz
inequality, provide a lower bound for the mean square integral of $E_2(T)$. 
The author proved in [7] the lower bound 
$$
\int_T^{2T}E^2_2(t)\d t \;\gg\;T^2,\eqno(1.12)
$$
which complements (1.3). However, in view of Theorem 1, it does
not appear likely that this procedure
can shed some new light on the behaviour of the
integral in (1.12).

Since the function $g(t)$ is $\Omega_\pm(1)$,  it means that it takes 
positive and negative values for some arbitrarily large values of $t$.
Thus it seems of interest to characterize the sets where $g(t)>0$
and $g(t)<0$. In this direction  we have ($\mu(\cdot)$ denotes measure)
the following result, which will be used in proving Theorem 2.

\bigskip
THEOREM 3. {\it There is a number } $\eta>0$ {\it such that $[T, 2T]$
contains a subset ${\Cal A}(T)$ in which $g(t) > \eta$ and 
$$
\mu({\Cal A}(T)) \;\gg_\eta\; T,\eqno(1.12)
$$
 and a subset ${\Cal B}(T)$ in which $g(t) < -\eta$ and}
$$
\mu({\Cal B}(T)) \;\gg_\eta\; T.\eqno(1.13)
$$

\bigskip
\head
2. Proof of Theorem 1
\endhead

Let us define
$$
G(T) := \int_0^T (E(t) - \pi)\d t.\eqno(2.1)
$$
If $P_4(x)$ is the polynomial appearing in the definition
of $E_2(T)$ and $Q_4(x) := P_4(x) + P_4'(x)$,
then integrating by parts we have 
$$
\eqalign{&
\int_T^{2T} E(t)E_2(t)\d t  = 
\int_T^{2T} (E(t)- \pi + \pi)E_2(t)\d t\cr&
= O(T^{3/2}) - \int_T^{2T}G(t)E_2'(t)\d t\cr&
= - \int_T^{2T}G(t)(|\zt|^4 - Q_4(\log t))\d t  + O(T^{3/2})\cr&
= - \int_T^{2T}G(t)|\zt|^4 \d t + O(T^{3/2}).\cr}\eqno(2.2)
$$
Proceeding with $E(t)$ in place of $E_2(t)$, we obtain (1.11).
Here we used the facts that 
$$ 
E_2(T) = O(T^{2/3}\log^8T),\quad
\int_0^T E_2(t)\d t = O(T^{3/2}),  \eqno(2.3)
$$
and
$$
G(T) = O(T^{3/4}),\quad\int_0^T G(t)\d t = O(T^{5/4}).\eqno(2.4)
$$ 
For a proof of the bounds in (2.3), see [4] or [12]. The 
bounds in (2.4) follow from  the explicit formula of
 Hafner--Ivi\'c [1], namely
$$
\eqalign{
G(t) &= 2^{-3/2}\sum_{n\le t}(-1)^nd(n)
n^{-1/2}\left({t\over2\pi n}+{1\over4}\right)^{-1/4}
\left(\roman {ar}\sinh\sqrt{\pi n\over 2t}\,\right)^{-2}\sin f(t,n)\cr&
-2\sum_{n\le c_0t}d(n)n^{-1/2}\left(\log{t\over2\pi n}\right)^{-2}
\sin \left(t\log{t\over2\pi n} - t - {\pi\over4}\right) + O(t^{1/4}),
\cr}\eqno(2.5)
$$
where $\;c_0 = 1/(2\pi) + 1/2 - \sqrt{1/4 + 1/(2\pi)}\;
= 0.019502\ldots\,$ and
$$
f(t,k) = 2t\,\roman {ar}\sinh\,\sqrt{\pi k\over2t} +
\sqrt{2\pi kt + \pi^2k^2} -{\textstyle{1\over4}}\pi,
 \;\roman {ar}\sinh\, x = \log(x + \sqrt{1+x^2}).
$$
By using the Cauchy--Schwarz inequality for integrals, 
the bound (e.g., see [2] for a proof)
$$
\int_0^T|\zt|^8\d t \ll T^{3/2}\log^{21/2}T   \eqno(2.6)
$$
and the mean theorem for Dirichlet polynomials (e.g., see [2, Chapter 5]), 
it is seen that the contribution of $\sum_{n\le c_0t}$ in (2.5) to the
right-hand side of (1.7) 
is $\ll T^{3/2}\log^{10}T$ (the exponent of the logarithm is not optimal, 
but it is unimportant). Simplifying the first sum
in (2.5) by Taylor's formula (truncating it at $n = \sqrt{T}$), 
we obtain from (2.2) the asympottic formula (1.7), namely
$$
I(T)  = 
- \int_T^{2T}t^{3/4}g(t)|\zt|^4 \d t + O(T^{3/2}\log^{10}T)\eqno(2.7)
$$
with $g(t)$ given by (1.6).
Since clearly 
$$
|g(t)| \le \hf \left({2\over\pi}\right)^{3/4}\sum_{n=1}^\infty d(n)n^{-5/4} 
=  \hf \left({2\over\pi}\right)^{3/4}\z^2({\textstyle{5\over4}}) = O(1),
$$
one obtains easily  from (2.7) and the weak bound
$$
\int_0^T|\zt|^4\d t \ll T\log^4T,\eqno(2.8)
$$
the upper bound
$$
I(T) \;=\; O(T^{7/4}\log^4T).\eqno(2.9)
$$
Although (2.9) improves (1.4), it is poorer than (1.5) of
Theorem 1, so that we must use different tools to obtain the
assertion of Theorem 1. To this end we appeal to the
following explicit formula of Ivi\'c--Motohashi (see [4] and [12]):
For $V^{1/2}\log^{-A}V \le \D \le V^{3/4}$ ($A>0$ is an arbitrary,
but fixed constant) 
$$
\eqalign{
\int_0^V I(T,\D)\d t &= VP_4(\log V) + O(\D\log^5V) \cr&
+ \pi\sqrt{\hf V}\sum_{j=1}^\infty \a_j\H c_j\cos
\bigl(\kk\log{\kk\over4{\roman e}V}\bigr)
{\roman e}^{-(\kk \D/2V)^2},\cr}\eqno(2.10)
$$
where $c_j \sim \k_j^{-3/2}$ as $\kk \to\infty$ and
$$
I(T,\D) \; := \;{1\over\D\sqrt{\pi}}\int_{-\infty}^\infty|\z(\hf + iT+iu)|^4
{\roman e}^{-(u/\D)^2}\d u.\eqno(2.11)
$$
For the definitions and properties of the spectral quantities
$\a_j, \kk$ and $H_j(\hf)$, see Y. Motohashi's monograph [13].
What will be needed here, besides (2.10) and (2.11), is essentially
the bound (cf. [13])
$$
\sum_{\k_j\le K}\a_j\H \ll K^2\log^3K. \eqno(2.12)
$$
Although (2.12) is not stated explicitly in [13], it follows when
one integrates the last formula on p. 130 (with $G = T^{3/4}$, say)
from $K$ to $2K$ with the help of the estimate for the sum in (2.12) in short
intervals, obtained recently by the author in [9]. An asymptotic formula for 
the sum in (2.12) has been obtained recently by  the author in [10]. This is
$$
\sum_{\k_j\le K}\a_j\H = K^2P_3(\log K) + R(K),\quad
R(K) = O(K^{5/4}\log^{25/2}K),  \eqno(2.13)
$$
where $P_3(x)$ is a cubic polynomial in $x$ with leading coefficient
equal to $4/(3\pi^2)$. 

\smallskip
The proof of (1.5) consists of three steps: the first is to show that  (2.2)
can be simplified to give (1.7). Then we show that $|\zt|^4$ can be
replaced by $I(t,\D)$ (with suitable $\D$) and permissible error. The
last step is to use the spectral decomposition (2.10) and obtain (1.5).
The bound in (1.5) is actually the limit of the method, set by
the condition $V^{1/2}\log^{-A}V \le \D \le V^{3/4}$ in (2.10). Namely
we wish $\D$ to be as small as possible, so any further improvements
of (1.5) will necessitate the  widening of this range, or obtaining 
another type of the asymptotic formula for the integral in (2.10).

\medskip
We proceed now with the proof. By using (2.6)
and the mean theorem for Dirichlet polynomials (e.g., see [2, Chapter 5]), 
it is seen that the contribution of $\sum_{n\le c_0t}$ in (2.5) to (2.2) 
is $\ll T^{3/2}\log^{10}T$ (the exponent of the logarithm is not optimal, 
but it is unimportant). Then we simplify the first sum
in (2.5) by Taylor's formula (truncating it at $n = \sqrt{T}$) to obtain
(1.7), as claimed.

\medskip
Before we go to the second step, let ($1 \ll N \ll \sqrt{T}$)
$$
g_N(t) := \hf \left({2\over\pi}\right)^{3/4}\sum_{n\le N}
(-1)^nd(n)n^{-5/4}\sin(\sqrt{8\pi nt} - {\textstyle{1\over4}}\pi).
\eqno(2.14)
$$
By using trivial estimation and (2.8) we have
$$\eqalign{&
\int_T^{2T}\sum_{n\ge N}
(-1)^nd(n)n^{-5/4}\sin(\sqrt{8\pi nt} - {\textstyle{1\over4}}\pi)
|\zt|^4\d t \cr&\ll TN^{-1/4}\log N\log^4T 
\ll T(\log T)^{7/2}\log\log T\cr} \eqno(2.15)
$$
by choosing 
$$
N \;=\; (\log T)^{2}.\eqno(2.16)
$$
Hence, in view of (1.6) and (2.7), it remains to prove that
$$
\int_T^{2T}g_N(t)|\zt|^4\d t \;\ll\; T(\log T)^{7/2}\log\log T.\eqno(2.17)
$$
Now we have, with $I(T,\D)$ given by (2.11),
$$
\eqalign{&
H(T,\D) := \int_T^{2T}g_N(t)(|\zt|^4 - I(t,\D))\d t\cr&
= {1\over\D\sqrt{\pi}}\int_{-\infty}^\infty
{\roman e}^{-(u/\D)^2}\int_T^{2T}g_N(t)(|\zt|^4 - |\z(\hf + it + iu)|^4)
\d t\,\d u\cr&
= {1\over\D\sqrt{\pi}}\int_{-\infty}^\infty
{\roman e}^{-(u/\D)^2}\int_T^{2T}g_N(t)\left(\int_u^0 
{\partial\over\partial x}
|\z(\hf + it + ix)|^4\d x\right)\d t\,\d u\cr&
= {1\over\D\sqrt{\pi}}\int_{-\infty}^\infty
{\roman e}^{-(u/\D)^2}\int_u^0\left(\int_T^{2T}g_N(t) 
{\partial\over\partial t}
|\z(\hf + it + ix)|^4\d t\right)\d x\,\d u\cr&
= {1\over\D\sqrt{\pi}}\int_{-\infty}^\infty
{\roman e}^{-(u/\D)^2}\int_{-|u|}^{|u|}O( 
|\z(\hf + iT + ix)|^4)\d x\,\d u\cr&
+ {1\over\D\sqrt{\pi}}\int_{-\infty}^\infty
{\roman e}^{-(u/\D)^2}\int_{-|u|}^{|u|}O( 
|\z(\hf + 2iT + ix)|^4)\d x\,\d u\cr&
-  {1\over\D\sqrt{\pi}}\int_{-\infty}^\infty
{\roman e}^{-(u/\D)^2}\int_u^0\int_T^{2T}g_N'(t) 
|\z(\hf + it + ix)|^4\d t\d x\,\d u.
\cr}
$$
The integrals over $u$ can be truncated at $|u| = \D\log T$ with a
negligible error. Since we have (see [4] and [12])
$$
\int_{T-G}^{T+G}|\zt|^4\d t \ll G\log^4T + TG^{-1/2}\log^CT\qquad(C > 0,
\, T^\e \le G \le T),
$$
and
$$
g_N'(t) \;\ll\; T^{-1/2}N^{1/4}\log N\qquad(T \le t \le 2T),\eqno(2.18)
$$
it follows that
$$
\eqalign{
H(T,\D) &\ll {1\over\D}\int\limits_{-\D\log T}^{\D\log T}
{\roman e}^{-(u/\D)^2}\d u
\Bigl(\int\limits_{T-\D\log T}^{T+\D\log T}|\zx|^4\d x\cr&
+ \int\limits_{2T-\D\log T}^{2T+\D\log T}|\zx|^4\d x\Bigr)\cr&
+ \D^{-1}T^{-1/2}N^{1/4}\log N
\int_{-\D\log T}^{\D\log T}{\roman e}^{-(u/\D)^2}|u|\d u\,T\log^4T\cr&
\ll \D\log^5T + T\D^{-1/2}\log^CT + T^{1/2}\D N^{1/4}\log N\log^4T\cr&
\ll T\log^{6-A}T \ll T(\log T)^{7/2}\log\log T \cr}\eqno(2.19)
$$
in view of (2.16), where we choose with $A > 0$ sufficiently large
(this is the lower bound in the
permissible range for which (2.10) holds)
$$
\D \;=\; T^{1/2}\log^{-A}T.
$$  

\medskip
We are now at the final step of the proof of Theorem 1. From (2.10) and
(2.11) we obtain, on integrating by parts (again $Q_4(x) = P_4(x)
+ P_4'(x)$),
$$
\eqalign{&
\int_T^{2T}g_N(t)I(t,\D)\d t =  \int_T^{2T}g_N(t)Q_4(\log t)\d t\cr&
+ O\left(\D\log^5T\int_T^{2T}|g_N'(t)|\d t\right) + O(\D\log^5T) +
O(T\D^{-1/2}\log^4T)\cr&
- \pi\int_T^{2T}g_N'(t)\left(\sqrt{\hf t}\sum_{\k_j\le T\D^{-1}\log T}
\a_j\H c_j\cos\bigl(\kk\log{\kk\over4{\roman e}t}\bigr){\roman e}^
{-(\kk\D/(2t))^2}\right)\d t,
\cr&}
$$
where we  used (2.12). We have
$$
\int_T^{2T}g_N(t)Q_4(\log t)\d t \ll T^{1/2}\log^4T
$$
by the first derivative test, and by using (2.18) we obtain
$$
\eqalign{&
\int_T^{2T}g_N(t)I(t,\D)\d t \cr&
\ll T^{1/2}\log^4T + T^{1/2}\D N^{1/4}\log N\log^5T + T\D^{-1/2}\log^4T
+ |K(T,\D)|,
\cr}\eqno(2.20)
$$
where we have set
$$\eqalign{
K(T,\D) &= \sum_{n\le N}(-1)^nd(n)n^{-3/4}\sum_{\k_j\le T\D^{-1}\log T}
\a_j\H c_j\times\cr&
\times \int_T^{2T}\cos(\sqrt{8\pi n t} - {\txt{1\over4}})
\cos\bigl(\kk\log{\kk\over4{\roman e}t}\bigr)
{\roman e}^{-(\kk\D/(2t))^2}\d t.\cr}
$$
We write the cosines as exponentials and note that the 
saddle point of the ensuing integral is at $t_0 = \kk^2/(2\pi n) \in [T,2T]$
for $\kk \asymp \sqrt{Tn}$. By the saddle point method (see e.g., [2])
the main contribution will be a multiple of ($F(t) = \sqrt{8\pi nt}
- \kk\log t$)
$$
\eqalign{&
\sum_{n\le N}(-1)^nd(n)n^{-3/4}\sum_{\kk\asymp\sqrt{Tn}}
\a_j\H c_j{\roman e}^{iF(t_0)}(F''(t_0))^{-1/2}
{\roman e}^{-(\kk\D/(2t_0))^2}\cr&
\ll \sum_{n\le N}d(n)n^{-7/4}\sum_{\kk\asymp\sqrt{Tn}}\a_j\H\cr&
\ll TN^{1/4}\log^3T\log N \ll T(\log T)^{7/2}\log\log T,\cr}\eqno(2.21)
$$
where we used again (2.12). Thus from (2.20) and (2.21)
we obtain  the bound in (2.17),
as asserted. This completes the proof of Theorem 1. In concluding, note
that the inner sum in (2.21) was estimated trivially. However, one
hopes that there is a lot of cancellation in such type of exponential
sum with $\a_j\H$. Indeed, it was conjectured by the author in [8]
that such a cancellation occurs, and it was heuristically justified
why one does expect this fact. Also there is hope to use the explicit 
expression which stands for the function $R(K)$ in the proof of the 
asympotic formula (2.13). The small improvement of the bound in
(1.5) of Theorem 1 over the bound in (2.9), which is relatively
not difficult to obtain, is precisely significant for this reason:
it does show that cancellation in a sum with $\a_j\H$ does occur.

\bigskip
\head
 3. Proof of Theorem 2 and Theorem 3
\endhead

\bigskip
We  shall first deal with Theorem 3, which is needed
for the proof of Theorem 2. The proof is based on the method used by the
author in [3]. Suppose $1 \ll H \ll T$. We
note that, by the first derivative test,
$$
\int_T^{T+H}g(t)\d t \;\ll\; \sqrt{T}  \eqno(3.1)
$$
holds uniformly in $H$, and proceed as follows. Let
$$ 
E = {1\over8}\left({2\over\pi}\right)^{3/2}
{\z^4({\textstyle{5\over2}})\over\z(5)}.
$$
Then we have ($EH$ comes from the terms $m = n$),  by the 
first derivative test,
$$
\eqalign{
&
\int_T^{T+H}g^2(t)\d t  = EH + O(\sqrt{T})\cr&
 + O\Bigl(\sum_{m,n=1;m\not=n}^\infty d(m)d(n)(mn)^{-5/4}\left|\int_T^{T+H}
\exp(i\sqrt{8\pi mt} \pm i\sqrt{8\pi nt})\d t\right|\Bigr)\cr&
= EH + O\Bigl(\sqrt{T}\sum_{m,n=1;m\not=n}^\infty d(m)d(n)(mn)^{-5/4}
|\sqrt{m}-\sqrt{n}|^{-1}\Bigr) + O(\sqrt{T})\cr&
= EH + O\Bigl(\sqrt{T}\sum_{r=1}^\infty \sum_{n=1}^\infty
(n(n+r))^{\e-5/4}r^{-1}(n+r)^{1/2}\Bigr)+ O(\sqrt{T})\cr&
= EH + O(\sqrt{T}).
\cr}\eqno(3.2)
$$ 
Thus (3.2) implies that
$$
\int_T^{T+H}g^2(t)\d t  \;\gg\;H, \eqno(3.3)
$$
provided that $H = D\sqrt{T}$ with a sufficiently large constant $D>0$.
We shall  show now that  there exists $\tau \in [T,T+H]$  
such that $g(\tau) > 2\eta$ for some constant $\eta > 0$ (and also a point
$\tau_1$ such that $g(\tau_1) < -2\eta$). Suppose on 
the contrary that $g(t) < \e$ if $g(t) > 0$ for any given $\e > 0$
(the case when  $g(t) > -\e$ if $g(t) < 0$ is treated analogously). Let
(cf. (1.8))
$$
g_+(t) = \max(g(t),\,0),\qquad g_-(t) = \min(g(t),\,0).\eqno(3.4)
$$
Since $g(t)$ is bounded we obtain, for some constants $C_1, C_2>0$, 
on using (3.1),
$$
\eqalign{
 \int_T^{T+H}g^2(t)\d t &=  \int_T^{T+H}g_+^2(t)\d t + 
\int_T^{T+H}(-g_-(t))^2\d t\cr&
\le H\e^2 + C_1\int_T^{T+H}(-g_-(t))\d t\cr&
= H\e^2 - C_1\int_T^{T+H}g(t)\d t + C_1\int_T^{T+H}g_+(t)\d t\cr&
\le H\e^2 +C_2\sqrt{T} + C_1\e H.\cr}
$$
However the above bound contradicts (3.3) if $\e$ is small enough and $D$
is large enough, since $H = D\sqrt{T}$. Hence there exists $\tau \in [T,T+H]$  
such that $g(\tau) > 2\eta$ for
some constant $\eta > 0$. Setting for brevity $c = \hf(2/\pi)^{3/4}$, using
$|\sin x - \sin y| \le |x - y| \;(x,y \in \RR)$, we have ($N \ge 3, x \gg 1$)
$$
\eqalign{&
|g(\tau + x) - g(\tau)| \cr&= \Big|c\sum_{n\le N}(-1)^nd(n)n^{-5/4}
\left(\sin(\sqrt{8\pi n(\tau + x)} - \pi/4) - \sin(\sqrt{8\pi n\tau} - \pi/4)
\right)\Big|
\cr&
+ O\left(\sum_{n>N}d(n)n^{-5/4}\right)\cr&
\ll \sum_{n\le N}d(n)n^{-5/4}n^{1/2}T^{-1/2}|x| + \sum_{n>N}d(n)n^{-5/4}\cr&
\ll T^{-1/2}|x|N^{1/4}\log N + N^{-1/4}\log N < \eta\cr}
$$
if $|x| \le L(T) :=  C\eta T^{1/2}N^{-1/4}\log N$ and 
$N = [\eta^{-8}]$, provided that $\eta$ is sufficiently small 
(which may be assumed), and $C>0$ is a suitable absolute constant. 
One obtains then
$$
|g(\tau  + x)| \ge g(\tau) - |g(\tau + x) - g(\tau)| > 
\eta, \quad |x| \le L(T).
$$
This means that every interval $[T_1, T_1 + D\sqrt{T}]\;(T \le T_1 \le 2T, D$
sufficiently large) contains a subinterval
of length $\gg_\eta \sqrt{T}$ in which $g(t) > \eta$. Consequently we
divide $[T,2T]$ into $\gg \sqrt{T}$ subintervals of the form
$[T + (j-1)D\sqrt{T}, T + jD\sqrt{T}],\;(j = 1,2,\ldots)$, and we obtain
that $g(t)>\eta$ on a set ${\Cal A}(T)$ satisfying
$\mu({\Cal A}(T)) \gg_\eta T$, as claimed.

\medskip

We pass now to the proof of Theorem 2. 
 The upper bounds in (1.9) and (1.10) follow easily
from the fact that $g(t)$ is bounded and that one has
$$
\int_T^{2T}|\zt|^{2k}\d t \;\ll\;T(\log T)^{k^2}\qquad(k = 1,2).
$$
For the lower bounds in (1.9) and (1.10) (they actually hold for 
any fixed $k\in\NN$) 
we use the well-known bound (see e.g., [2] or [4])
$$
\int_T^{T+C\sqrt{T}}|\zt|^{2k}\d t \;\gg\;\sqrt{T}(\log T)^{k^2}
\qquad(C > 0,\,k \in \NN).
\eqno(3.5)
$$
In each interval of length $\ge C\sqrt{T}$ in the proof of Theorem 2
(where we had $g(t) > \eta$) we use (3.5) and gather the resulting 
lower bounds to obtain the lower bound in (1.9). 
The proof of the lower bound in (1.10) is analogous.




\Refs
\bigskip
\item{[1]} J. L. H a f n e r and A. I v i \'c, {\it
On the mean square of the Riemann zeta-function on the critical line}, 
J. Number Theory   {\bf 32}(1989), 151-191.

\smallskip
\item{[2]} A. I v i \'c, {\it  The Riemann zeta-function}, 
John Wiley \& Sons, New York, 1985.

\smallskip
\item {[3]}  A. I v i \'c, {\it Large values of certain number-theoretic
error terms}, Acta Arith. {\bf56}(1990), 135-159.

\smallskip
\item {[4]}  A. I v i \'c, {\it  Mean values of the Riemann zeta-function},
LN's {\bf82},  Tata Institute of Fundamental Research, Bombay, 1991 
(distr. by Springer Verlag, Berlin etc.).

\smallskip
\item{[5]}  A. I v i \'c, {\it On the fourth moment of the Riemann 
zeta-function}, 
Publications Inst. Math. (Belgrade) {\bf57(71)} (1995), 101-110.

\smallskip
\item{[6]} A.  I v i \'c, {\it The Mellin transform and the Riemann 
zeta-function}, 
Proceedings of the  Conference on Elementary and Analytic Number Theory 
(Vienna, July 18-20, 1996), Universit\"at Wien \& Universit\"at f\"ur
Bodenkultur, Eds. W.G. Nowak and J. Schoi{\ss}engeier, Vienna 1996,
pp. 112-127.

\smallskip
\item{[7]} A. I v i \'c, {\it On the error term for the fourth moment of 
the Riemann zeta-function}, J. London Math. Soc., 
(2){\bf60}(1999), 21-32.

\smallskip
\item{[8]} A. I v i \'c, {\it
On some conjectures and results for the Riemann zeta-function
and Hecke series}, Acta Arith. {\bf109}(2001), 115-145.

\smallskip
\item{ [9]} A. I v i \'c, 
{\it On sums of Hecke series in short intervals}, J. Th\'eorie des
Nombres Bordeaux {\bf14}(2001), 554-568.

\smallskip
\item{ [10]} A. I v i \'c, {\it On the moments of Hecke series
at central points}, Functiones et Approximatio {\bf30}(2002), 49-82.

\smallskip
\item{ [11]} A. I v i \'c and Y. M o t o h a s h i, {\it The mean 
square of the error term for the fourth moment of the zeta-function}, 
 Proc. London Math. Soc. (3){\bf66}(1994), 309-329.

\smallskip

\item {[12]}  A. I v i \'c and Y. M o t o h a s h i, {\it  The fourth 
moment of the Riemann zeta-function}, J. Number Theory {\bf 51}(1995), 16-45.

\smallskip
\item {[13]} Y. M o t o h a s h i, {\it Spectral theory of the Riemann
zeta-function},  Cambridge University Press, Cambridge, 1997.

\bigskip\bigskip

Aleksandar Ivi\'c

Katedra Matematike RGF-a

Universiteta u Beogradu

\DJ u\v sina 7, 11000 Beograd, Serbia

{\tt aivic\@rgf.bg.ac.yu}

\endRefs
\bye